\documentclass[leqno,11pt, english]{amsart}
\usepackage[usenames,dvipsnames]{color}
\usepackage[colorlinks=true,linkcolor=Red,citecolor=Green]{hyperref}
\usepackage{amsmath,amssymb,amsthm,graphicx,url}
\usepackage{epstopdf}
\usepackage{color}
\usepackage{dsfont}
\usepackage[utf8]{inputenc} 
\usepackage[T1]{fontenc}
\newcommand{\xqedhere}[1]{%
    \rlap{%
         \hbox to#1{%
           \hfil
           \llap{%
               \ensuremath{\square}
           }%
       }%
   }%
}

\def\pasdegrille{\let\grille = \pasgrille}

\def\aat#1#2#3{
\divide \dimen1 by 48 \dimen3=\dimen1 \multiply \dimen1 by #1
\advance \dimen1 by -\dimen3 \divide \dimen1 by 101 \multiply
\dimen1 by 100 \divide \dimen2 by \count11 \multiply \dimen2 by #2
\setbox0=\hbox{#3}\ht0=0pt\dp0=0pt
  \rlap{\kern\dimen1 \vbox to0pt{\kern-\dimen2\box0\vss}}\dimen1= \wd1
\dimen2=\ht1}
\def\pasgrille{
\count12= \dimen1 \divide \count12 by 50 \divide \dimen2 by
\count12 \count11 =\dimen2 \ \divide \dimen1 by 48
\setlength{\unitlength}{\dimen1} \smash{\rlap{\ }} \dimen1= \wd1
\dimen2=\ht1 }
\def\grille{
\count12= \dimen1 \divide \count12 by 50 \divide \dimen2 by
\count12 \count11 =\dimen2 \ \divide \dimen1 by 48
\setlength{\unitlength}{\dimen1}
\smash{\rlap{\graphpaper[1](0,0)(50, \count11)}} \dimen1= \wd1
\dimen2=\ht1 }

\pasdegrille

\newcommand{\be}{\begin{equation}}
\newcommand{\ee}{\end{equation}}

\theoremstyle{plain}

\newtheorem{thm}{Theorem}

\newtheorem{prop}{Proposition}[section]

\theoremstyle{definition}

\newtheorem{rem}[prop]{Remark}
\newtheorem{defn}[prop]{Definition} 

\numberwithin{equation}{section}

\def\squarebox#1{\hbox to #1{\hfill\vbox to #1{\vfill}}}

\usepackage{amsxtra}

\ifx\pdfoutput\undefined
  \DeclareGraphicsExtensions{.pstex, .eps}
\else
  \ifx\pdfoutput\relax
    \DeclareGraphicsExtensions{.pstex, .eps}
  \else
    \ifnum\pdfoutput>0
      \DeclareGraphicsExtensions{.pdf}
    \else
      \DeclareGraphicsExtensions{.pstex, .eps}
    \fi
  \fi
\fi

\title[Simutaneous control]{Simultaneous control for the heat equation with Dirichlet and Neumann boundary conditions}

\author[N. Burq]{Nicolas Burq}
\address{Universit{\'e} Paris-Saclay, Math{\'e}matiques, UMR 8628 du CNRS, B{\^a}t 307, 91405  Orsay Cedex, France,   and Institut Universitaire de France}
\email{Nicolas.burq@universite-paris-saclay.fr}
\author[I. Moyano]{Ivan Moyano}
\address{Universit\'e de Nice Sophia-Antipolis
Parc Valrose, Laboratoire J.A. Dieudonn\'e,
UMR 7351 du CNRS
06108 NICE Cedex 02
FRANCE}
\email{Ivan.Moyano@unice.fr}

\usepackage{amssymb}
\usepackage{amsmath, amsthm, amsopn, amsfonts}

\def\11{{\rm 1~\hspace{-1.4ex}l} }

\begin{document}

\begin{abstract}
It is well known that both the heat equation with Dirichlet or Neumann boundary conditions are null controlable as soon as the control acts in a non trivial domain (i.e. a set of positive measure, see~\cite{LeRo95-1, mi10, PhWa13, ApEs13, BM}. In this article, we show that for any couple of initial data $(u_0, v_0)$ we can achieve the null control for both equations (Dirichlet and Neumann boundary conditions respectively)  simultaneously  {\em with the same control function for both equations}. \begin{center}
 
  \rule{0.6\textwidth}{.4pt} \end{center}\end{abstract}   

\ \vskip -1cm \noindent\hfil\rule{0.9\textwidth}{.4pt}\hfil \vskip 1cm 
 \maketitle   
 \section{Introduction and main results}  
Let us consider a smooth bounded domain $\Omega \subset \mathbb{R}^d$ and let $\omega \subset \Omega$ be a subset of positive measure $|\omega| >0$ and the following internal simultaneous controlability problem
\begin{equation}\label{heat}
\left \{ \begin{aligned} (\partial_t - \Delta ) u & = f 1_{(0,T) \times \omega},  \qquad u \mid_{\partial \Omega} = 0, \qquad  u \mid_{t=0} = u_0, \\
(\partial_t - \Delta ) v & = f 1_{(0,T)\times \omega},  \qquad \partial_{\nu} v \mid_{\partial \Omega} = 0, \qquad  v \mid_{t=0} = v_0. 
\end{aligned}
\right.
\end{equation}

\begin{defn}
We shall say that the heat equation in $\Omega$ is {\em simultaneously null controlable} with Dirichlet and Neumann boundary conditions if 
for any $(u_0, v_0) \in L^2(\Omega)$ there exists $f \in L^2((0,T) \times \omega)$ such that the solution of the system~\eqref{heat} satisfies 
$$ u \mid_{t>T} =0, \quad v\mid_{t>T} =0.$$
\end{defn}

The question of simultaneous controllability of various partial differential equations has been raised in the literature (see for example \cite{AvdoninPandolfi,SimultaneousSchrodinger,Ammar,FdezCara,Benabdallah} and the recent work \cite{ArarunaSimultaneous} for more references on the subject), especially when the system involves some transmission mechanisms between the equations allowing to reduce the number of commands. The interest of our problem lies on the fact that both heat equations in (\ref{heat}) exhibit relatively independent dynamics and yet they can be steered to zero using \emph{exactly the same control}.  As no coupling exists between these two equations, this simultaneous controlability is at first glance counter intuitive. Yet, by considering the two new unknowns $w_1= u+v, w_2 = u-v$, the simultaneous controlability reduces to the controlability of $(w_1, w_2)$, with a control acting only on the $w_1$ component of the system. Notice that the system $(w_1, w_2)$ is now coupled at the boundary by the transmission conditions
$$ \partial_\nu w_1 \mid_{\partial \Omega} =\partial_\nu w_2 \mid_{\partial \Omega} , \qquad w_1 \mid_{\partial \Omega}= -w_2 \mid_{\partial \Omega} ,$$ 
 and we need to show that this coupling is sufficient. However, our strategy will follow a more direct path (the double manifold) and will not study {\em per se} this transmission problem.

\subsection{Simultaneous controllability}

Our first result is the following
\begin{thm} \label{thheat}
Let $T>0, \omega \subset \Omega$ of positive measure $|\omega| >0$. Then the heat equation in $\Omega$ is {simultaneously null controlable} with Dirichlet and Neumann boundary conditions.
\end{thm}

\begin{rem} It is classical that both the heat equation with Dirichlet or with Neumann conditions are null controlable. The novelty in Theorem~\ref{thheat} lies precisely on the fact that the null controlability can be achieved for any initial data $(u_0,  v_0) \in L^2(\Omega)\times L^2(\Omega)$  {\em with the same control} for both. 
\end{rem}

\begin{rem}
The proof we give below relies on the doubling manifold approach from~\cite{BM}. This approach is very robust and allows rough domains (of class $W^{2, \infty}$, for instance) and rough space-dependent Laplace operators (Lipschitz coefficients). 
$$ \Delta = \frac 1 {\kappa(x)} \sum_{i,j} \partial_{x_i} g^{i,j}(x) \kappa (x) \partial_{x_j},$$
where we assume that the coefficients $\kappa, g$ are Lipschitz  and that $g$ is uniformly elliptic.
As will appear clearly,  the proof (which is very simple once the results in~\cite{BM} were established) shows that all the control results from~\cite{BM} are true with the same control functions for the Dirichlet and Neumann heat equations.
\end{rem}

\begin{rem}
It is an interesting question whether similar results might hold for the wave equation. We plan to address this question in a forthcoming paper. However, in this case, the analysis is much more involved and we do not expect to get such a general answer.
\end{rem}

\subsection{Simultaneous controllability and spectral inequalities}

Let $(e^D_{\lambda})_{\lambda}$ be the spectral family associated to the Laplacian in $\Omega$ with Dirichlet conditions, i.e.,
\begin{equation}
-\Delta e^D_{\lambda} = \lambda^2 e^D_{\lambda}, \qquad e^D_{\lambda}\vert_{\partial \Omega} = 0  
\label{eq:SpectralDirichlet}
\end{equation} and let $(e^N_{\mu})_{\mu}$ be the spectral family associated to the Laplacian in $\Omega$ with Neumann conditions,
\begin{equation}
-\Delta e^N_{\mu} = \mu^2 e^N_{\mu}, \qquad \partial_{\nu} e^N_{\mu}\vert_{\partial \Omega} = 0. 
\label{eq:SpectralNeumann} 
\end{equation} Let $\omega \subset \Omega$ be a non-empty subset. In the spirit of \cite{JeLe96}, one can show that the spectral families $(e^D_{\lambda})_{\lambda}$ and $(e^N_{\mu})_{\mu}$ enjoy a concentration property on the subsets $\omega$ as long as they are not too small. In \cite{BM} the authors shows the following: given $\omega$ with $|\omega |>0$, there exist constants $C,D$ such that all the spectral truncations
\begin{equation*}
\Pi^D_{\Lambda} u := \sum_{\lambda \leq \Lambda} u_{\lambda} e^D_{\lambda}, \qquad u \in L^2(\Omega), \quad \Lambda > 0, 
\end{equation*} satisfy the estimate (cf. \cite[Theorem 1]{BM}): 
\begin{equation}\label{BORNE}
 \| \Pi^D_{\Lambda} u \|_{L^\infty(\Omega)} \leq C e^{C \Lambda} \| \mathds{1}_{\omega} \Pi^D_{\Lambda} u \|_{L^1(\Omega)}, \qquad \forall u \in L^2(\Omega).
\end{equation} The analogous estimate also holds for the spectral truncations of $(e^N_{\mu})_{\mu}$, defined by
\begin{equation*}
\Pi^N_{\Lambda} u := \sum_{\mu \leq \Lambda} u_{\mu} e^N_{\mu}, \qquad u \in L^2(\Omega), \quad \Lambda > 0. 
\end{equation*} In this note we show that these spectral inequalities also hold \textit{simultaneously} (i.e. we can estimate each spectral projector by the sum on arbitrary small set of positive measures).

\begin{thm}\label{spectral}
Let $\omega \subset \Omega$ with $ |\omega| \geq m$. There exist $C, D>0$ such that for any $\Lambda >0$,   
we have  
\begin{equation}\label{BORNE SIMULTANEE}
 \| \Pi^D_{\Lambda} u \|_{L^\infty( \Omega)}+ \|\Pi^N_{\Lambda} v \|_{L^\infty(\Omega)} \leq C e^{C \Lambda} \| \mathds{1}_{\omega}(\Pi^D_{\Lambda} u + \Pi^N_{\Lambda} v) \|_{L^1(\omega)}, \qquad \forall u,v \in L^2(\Omega).
\end{equation}

\end{thm}

\section{Double manifold and spectral estimates}
\label{sec:double mfld}

In this section we recall a result from \cite{BM}, which allows to  glue any given manifold $M$ with a copy of itself along its boundary, in order to produce a \textit{double manifold} without boundary. This will be a crucial point in the  analysis below. \par

\subsection{The double manifold}

Let $(M,g)$ be a compact Riemannian manifold of class $C^1 \cap W^{1,\infty}$.

Let $\Delta$  be the  Laplace-Beltrami operator on $M$ and let $(e_k)$ be a family of eigenfunctions of $- \Delta$, with eigenvalues $\lambda_k^2 \rightarrow + \infty$ forming a Hilbert basis of $L^2(M)$. 
$$ -\Delta e^{D;N}_k = \lambda_k ^2 e_k, \qquad e^D_k \mid_{\partial M} =0 \text{ (Dirichlet condition) or } \partial_\nu e^N_k \mid_{\partial M} =0 \text{ (Neumann condition)}.$$ Let be $\widetilde{M}$ the double space made of two copies of $\overline{M}$
\begin{equation*}
\widetilde{M} = \overline{M} \times \{-1,1\}/\partial M,
\end{equation*} where we identified the points on the boundary, $(x,-1)$ and $(x,1)$, $x\in \partial M$. In the double manifold $\widetilde{M}$ we have the following result.

\begin{thm}[The double manifold, \protect{\cite[Theorem 7]{BM}}]\label{double}
Let $g$ be given. There exists a $W^{2, \infty}$ structure on the double manifold $\widetilde{M}$, a metric $\widetilde{g}$  of class $W^{1, \infty}$ on $\widetilde{M}$, and  a density $\widetilde{\kappa}$ of class $W^{1, \infty}$ on $\widetilde{M}$ such that the following holds.
\begin{itemize}
\item The maps
$$ i^\pm  x\in M  \rightarrow (x, \pm 1) \in \widetilde{M} = M \times \{\pm 1\} / \partial M$$ are isometric embeddings.
\item The density induced on each copy of $M$ is the density $\kappa$,
$$\widetilde{\kappa} \mid_{M\times\{\pm1 \}} = \kappa.
$$
\item For any eigenfunction $e$ with eigenvalue $\lambda^2$ of the Laplace operator $-\Delta = - \frac 1 {\kappa} \text{div } g^{-1} \kappa \nabla$ with Dirichlet or Neumann boundary conditions, there exists an eigenfunction $\widetilde{e}$ with the same eigenvalue $\lambda$ of the Laplace operator $-\Delta = - \frac 1 {\widetilde{\kappa} }\text{div } \widetilde{g}^{-1} \widetilde{\kappa} \nabla$ on $\widetilde{M}$ such that 
\begin{equation}\label{ext}
\widetilde{e} \mid_{M \times \{1\}} = e, \quad \widetilde{e} \mid_{M \times \{-1\}} = \begin{cases}  -e \quad &(\text{Dirichlet boundary conditions}),\\
 e \qquad &(\text{Neumann boundary conditions}). \end{cases}
\end{equation}
\item Conversely, there exists a Hilbert basis of $L^2(\widetilde M)$ composed of eigenfunctions of the Laplace operator $\widetilde{\Delta}$ which are either odd extensions of Dirichlet Laplace eigenfunctions in $M$ or even extensions of Neumann Laplace eigenfunctions in $M$.
\end{itemize}
\end{thm}
\begin{rem} The last property was not stated explicitely in~ \cite[Theorem 7]{BM}, but it is straightforward as the vector space generated by such eigenfunctions is clearly dense in $L^2( \widetilde{M})$.
\end{rem}

\subsection{Spectral projector on the double manifold and proof of Theorem~\ref{spectral}}
Let us denote by $\widetilde{\Pi}_\Lambda$ the spectral projector on the manifold $\widetilde{M}$. Let $u, v \in L^2(M)$ and define the function 
\begin{equation}\label{extension}
\widetilde{u} (x,1)= (u+v) (x), \qquad \widetilde{u} (x,-1)= (-u+v) (x).
\end{equation}
Clearly if 
$$ u = \sum_k u_k e_k^D, \qquad v = \sum_k v_k e^N_k,$$
we get
$$ \widetilde{u} = \sum_k u_k \widetilde{e}^D_k + v_k \widetilde{e}^N_k.$$
 According to the reflection principle of the previous section, we can link 
 the Dirichlet and Neumann spectral projectors on $M$ and the spectral projector on $\widetilde{M}$ by the relation 
\begin{equation}\label{link} \widetilde{\Pi}_\Lambda (\widetilde{u}) (\cdot, 1) = \Pi_\Lambda^D (u(\cdot)) + \Pi^N_\Lambda (v(\cdot)), \qquad \widetilde{\Pi}_\Lambda (\widetilde{u} )(\cdot, -1) = -\Pi_\Lambda^D (u(\cdot)) + \Pi^N_\Lambda (v(\cdot)).
\end{equation}

\begin{thm}[\protect{\cite[Theorem 1]{BM}}] \label{spectralbis}
Let $\widetilde{\omega} \subset \widetilde{M}$ with positive Lebesgue measure.
Then, there exists $C>0$ such that for anu $\Lambda >0$ and any $\widetilde{u}\in L^2( \widetilde{M})$, we have 
\begin{equation}\label{spec-proj}
\| \widetilde{\Pi} _\Lambda \widetilde{u}\|_{L^{\infty} (\widetilde{M})} \leq  C e^{C\Lambda} \| 1 _{\widetilde{\omega}} \widetilde{\Pi} _\Lambda \widetilde{u} \|_{L^1(\widetilde{\omega})}.
\end{equation}
\end{thm}
We can now prove Theorem~\ref{spectral}. Indeed, let $\omega\subset M$ of positive Lebesgue measure. Let $\widetilde{\omega} = \omega \times \{1\}$.  According to Theorem~\ref{spectralbis} and~\eqref{link}, we get for any $u, v\in L^2(M)$, 
\begin{multline}
\| \Pi ^D_\Lambda u \|^2_{L^\infty(M)} + \| \Pi ^N_\Lambda v \|^2_{L^\infty(M)} = \| \Pi ^D_\Lambda u +  \Pi ^N_\Lambda v\|^2_{L^\infty(M)} + \| \Pi ^D_\Lambda u -\Pi ^N_\Lambda v\|^2_{L^\infty(M)}\\
= \| \widetilde{\Pi} _\Lambda \widetilde{u} \|^2_{L^\infty(\widetilde{M})} \leq C e^{C\Lambda} \| 1_{\widetilde{\omega} }\widetilde{\Pi} _\Lambda \widetilde{u} \|^2_{L^1(\widetilde{\omega})}\\
= C e^{C\Lambda} \| 1_{\omega\times\{1\}} \widetilde{\Pi} _\Lambda \widetilde{u} \|^2_{L^1(\widetilde{\omega})}= C e^{C\Lambda} \| 1_{\omega}\Pi ^D_\Lambda u +  \Pi ^N_\Lambda v \|^2_{L^1({\omega})}.
\end{multline}
\subsection{Control and the double manifold}
To prove our control result, we could just apply the spectral projector estimate we just proved and some functional analysis. Here we prefered to prove the result directly on the double manifold. We start with 
\begin{thm}[\protect{\cite[Theorem 2] {BM}}]\label{heat-double}
Let $\widetilde{\omega} \subset \widetilde{M}$ be a measurable set with $\vert \widetilde{\omega} \vert >0$. Then, for every $T>0$ and every $\tilde{u}_0 \in L^2(M)$, there exists $\tilde{f} \in L^2((0,T)\times \tilde{\omega})$ such that the solution to the heat equation on $\widetilde{M}$ satisfies 
\begin{equation*}
\tilde{u}\vert_{t \geq T} = 0.
\end{equation*}
\label{thm:control double}
\end{thm}
We can now prove Theorem \ref{thheat}. 
For any $u, v \in L^2(M)$, let us define $\widetilde{u}$ by~\eqref{extension}, and for any $\omega \subset M$ of positive measure, let $\widetilde{\omega} = \omega \times \{1\}$. 
 According to Theorem \ref{heat-double}, for every $T>0$, there exists $\tilde{f} \in L^2((0,T)\times \widetilde{\omega})$ such that 
\begin{equation}
(\partial_t - \tilde{\Delta} ) \widetilde{U}  = \tilde{f} \mathds{1}_{(0,T) \times \omega},   \qquad  \widetilde{U} \mid_{t=0} = \tilde{u}, \qquad  \widetilde{U}\vert_{t \geq T} = 0.
\label{eq:heat doubled}
\end{equation} Let us define next 
\begin{equation*}
u(t,x) = \widetilde{U}(t,x,1) - \widetilde{U}(t,x,-1), \qquad  v(t,x) = \widetilde{U}(t,x,1) + \widetilde{U}(t,x,-1),
\end{equation*} where $\widetilde{U}$ is defined by (\ref{eq:heat doubled}). Notice that $u$ clearly satisfies the Dirichlet boundary condition while $v$ satisfies the Neumann boundary condition.  This second condition is not obvious but comes from the construction of the double manifold in~\cite{BM}. Indeed, in our construction, we defined  normal coordinate system near any point in the boundary of $M$ such that $M = \{ x_n >0\}$, and then we glued the two copies defined by $M\times \{1\}  = \{ x_n >0\} $, $M\times\{-1\} = \{ x_n <0\} $ by the relation  
$$ (x,1) = (x_n, x',1), (x,-1) = (-x_n, x', -1), $$
which implies 
$$\partial _\nu v = \partial_{x_n} (  \widetilde{U}(t,x,1) + \widetilde{U}(t,x,-1))\mid_{x_n=0} =  \partial_{x_n} (  \widetilde{U})(t,x,1)) - \partial_{x_n} (\widetilde{U})(t,x,-1))  =0.
$$   
Now, by definition of $u$ and $v$ we have 
\begin{equation*}
(\partial_t - \Delta ) u  = \tilde{f}(t,x,1) \mathds{1}_{(0,T) \times \omega} - \tilde{f}(t,x,-1) \mathds{1}_{(0,T) \times \tilde{\omega}} = f(t,x) \mathds{1}_{(0,T) \times \omega}, 
\end{equation*} as $\tilde{f}(t,x,-1) \mathds{1}_{(0,T) \times \tilde{\omega}} = 0$ by the choice of $\widetilde{\omega} = \omega \times \{1\}$. By the same token, we have 
\begin{equation*}
(\partial_t - \Delta ) v  = \tilde{f}(t,x,1) \mathds{1}_{(0,T) \times \omega} + \tilde{f}(t,x,-1) \mathds{1}_{(0,T) \times \tilde{\Omega}} = f(t,x) \mathds{1}_{(0,T) \times \omega}. 
\end{equation*} As a consequence, $u$ and $v$ solve (\ref{heat}) with control $f \mathds{1}_{(0,T) \times \omega}$. Finally, using (\ref{eq:heat doubled}), we get
\begin{equation*}
u \vert_{t \geq T} = 0, \qquad v \vert_{t \geq T} = 0, 
\end{equation*} which ends the proof.

\bibliographystyle{plain}

\end{document}